\documentclass[11pt]{amsart}
\usepackage{latexsym}
\usepackage{amssymb,amscd,amsthm,amsfonts,verbatim,amsopn,color}
\usepackage[dvips]{epsfig}
\def\Z{{\mathbb Z}}

\def\R{{\mathbb R}}
\def\N{{\mathbb N}}

\def\P{{\mathbb P}}

\def\AA{{\mathcal A}}

\def\LL{{\mathcal L}}
\def\MM{{\mathcal M}}

\def\Ss{{\mathcal S}}

\newcommand{\lra}{\longrightarrow}
\newcommand{\stab}{{\rm stab}}

\newcommand{\fix}{\text{\rm Fix}}
\newcommand{\ra}{\rightarrow}

\newcommand{\car}{\circlearrowright}

\newtheorem{thm}{Theorem}[section]
\newtheorem{df} [thm]{Definition}
\newtheorem{lm}  [thm]{Lemma}

\newtheorem{prop}[thm]{Proposition}

\newtheorem{rem}[thm]{Remark}

\newtheorem{expl}[thm]{Example}
\numberwithin{equation}{section}

\newenvironment{pf}{\noindent {\bf Proof.}}{\hfill $\Box$\vspace{0.3cm}}
\newenvironment{explrm}{\begin{expl} \rm}{\end{expl}}
\newenvironment{remrm}{\begin{rem} \rm}{\end{rem}}

\begin{document}

\title[A desingularization of real differentiable actions of finite groups]
{A desingularization of real \\ differentiable actions of finite groups}

\author{Eva Maria Feichtner \mbox{ }\& \mbox{ }Dmitry N.\ Kozlov}

\address{
Department of Mathematics, ETH Zurich, 8092 Zurich, Switzerland}
\email{feichtne@math.ethz.ch}

\address{
Department of Mathematics, KTH Stockholm, 100 44 Stockholm, Sweden}
\curraddr{
Institute for Theoretical Computer Science, ETH Zurich, 8092 Zurich, 
Switzerland}
\email{kozlov@math.kth.se, dkozlov@inf.ethz.ch}


\begin{abstract}
  We provide abelianizations of differentiable actions of finite
  groups on smooth real manifolds.  De~Concini-Procesi wonderful
  models for (local) subspace arrangements and a careful analysis of
  linear actions on real vector spaces are at the core of our
  construction.  In fact, we show that our abelianizations have
  stabilizers isomorphic to elementary abelian $2$-groups, a setting
  for which we suggest the term {\em digitalization\/}. As our main
  examples, we discuss the resulting digitalizations of the
  permutation actions of the symmetric group on real linear and
  projective spaces.
\end{abstract}

\maketitle

\section{Introduction}

Abelianizations of finite group actions on complex manifolds appeared
prominently in the work of Batyrev~\cite{Ba}, and a connection to the
wonderful arrangement models of De~Concini and Procesi was observed by
Borisov and Gunnells~\cite{BG}.  The authors of the present article
have previously presented a detailed study of the key example over the
reals, the abelianization of the permutation action of the symmetric
group $\Ss_n$ on $\R^n$ given by the maximal De~Concini-Procesi model
of the braid arrangement, cf.~\cite{FK2}. In particular, it was shown
that stabilizers of points on the arrangement model are elementary
abelian~$2$-groups. We suggest to call an abelianization with this
property a {\em digitalization\/} of the given action.

In the present article, we extend our analysis in two steps.  First,
for any linear action of a finite group on a real vector space, we
define an arrangement of linear subspaces whose maximal
De~Concini-Procesi model we then show to be a digitalization of the
given action.  Second, we proceed by analysing differentiable actions
of finite groups on smooth real manifolds. We propose a locally finite
stratification of the manifold by smooth submanifolds and, observing
that this stratification is actually a local subspace arrangement, we
show that the associated maximal De~Concini-Procesi model is a
digitalization of the given action.

We present examples in the linear and in the non-linear case. First,
we consider the permutation action of the symmetric group~$\Ss_n$ on
$\R^n$, and we find that our arrangement construction yields the rank
$2$~truncation of the braid arrangement. The resulting digitalization
is the one discussed in~\cite{FK2}. It would be interesting to extend
this result to the action of an arbitrary reflection group. As a
non-linear example, we consider the action of $\Ss_n$ on $\R\P^{n-1}$
given by projectivizing the real permutation action on~$\R^n$. We show
that our manifold stratification, in this case, coincides with the
rank $2$~truncation of the projectivized braid arrangement. The
resulting digitalization thus is the maximal projective
De~Concini-Procesi model for the braid arrangement.

We give a short overview on the material presented in this article: In
Section~\ref{sec_review} we provide a review of De~Concini-Procesi
arrangement models in an attempt to keep this exposition fairly
self-contained. Our main results are presented in
Section~\ref{sec_digital}. In~\ref{ssec_lin} we propose a
digitalization for any given linear action of a finite group on a real
vector space; in~\ref{ssec_mnfds} we extend our setting to
differentiable actions of finite groups on smooth real manifolds.
Section~\ref{sec_appl} is focused on examples.  We work out the
details of the proposed digitalizations for the real permutation
action in~\ref{ssec_MaxAn}, and for the permutation action on real
projective spaces in~\ref{ssec_MaxPAn}.


\section{A review of De~Concini-Procesi arrangement models}
\label{sec_review}

\subsection{Arrangement models} \label{ssec_arrgtmodels}

We review the construction of De~Concini-Procesi arrangement models as 
presented in~\cite{DP1}. Moreover, we recall an encoding of points in
maximal arrangement models from~\cite{FK2} that will be  crucial for the 
technical  handling of stabilizers (cf.\ \ref{ssec_Garrgtmodels}).  

\subsubsection{The model construction} \label{sssec_modelconstr}
Let $\AA$ be a finite family of linear subspaces in some real or 
complex vector space~$V$. The
combinatorial data of such subspace arrangement is customarily recorded by
its {\em intersection lattice\/} $\LL\,{=}\,\LL(\AA)$, the partially ordered
set of intersections among subspaces in $\AA$ ordered by reversed
inclusion. We agree on the empty intersection to be the full space
$V$, represented by the minimal element $\hat 0$ in the lattice. 
We will frequently use $\LL_{>\hat 0}$ to denote $\LL\,{\setminus}\,\{\hat 0\}$.

There is a family of arrangement models each coming from the choice of
a certain subset of the intersection lattice, so-called {\em building sets\/}.
For the moment we restrict our attention to the maximal model among those, 
which results from choosing the whole intersection lattice as a building set.

We give two alternative descriptions for the maximal
De~Concini-Procesi model of $\AA$.    
Consider the following map on the complement $\MM(\AA)\,{:=}\,
V\,{\setminus}\,\bigcup \AA$ of the arrangement,
\begin{equation} \label{eq_modelconstr}
    \Psi: \quad \MM(\AA) \, \, \longrightarrow \, \, 
                    V\, \times 
                    \, \prod_{X\in\LL_{>\hat 0}}\, \P(V/X)\, ,
\end{equation}
where $\Psi$ is the natural inclusion into the first factor and
the natural projection to the other factors restricted to $\MM(\AA)$.
Formally,
\[
   \Psi (x)\, \, = \, \, (\,x\,, (\langle x,X \rangle/X)_{X\in \LL_{> \hat 0}})\, ,
\]
where $\langle \cdot,\cdot \rangle$ denotes the linear span of
subspaces or vectors, respectively, and $\langle x,X \rangle/X$ is interpreted 
as a point in $\P(V/X)$ for any $X\,{\in}\,\LL_{> \hat 0}$.  

The map $\Psi$ defines an embedding of $\MM(\AA)$ into the product on the 
right hand side of~(\ref{eq_modelconstr}). 
The closure of its image, $Y_{\AA}:= \overline{{\rm im}\Psi}$, 
is the {\em maximal De Concini-Procesi model\/} of the arrangement~$\AA$.
If we want to stress the ambient space of the original arrangement, we will
use the notation $Y_{V,\AA}$ for $Y_{\AA}$.

Alternatively, one can describe $Y_{\AA}$ as the result of successive
blowups of strata in $V$. Consider the stratification of $V$ given by the
linear subspaces in $\AA$ and their intersections. Choose some linear 
extension of the opposite order in $\LL$.
Then, $Y_{\AA}$ is the result of successive blowups of strata, respectively
proper transforms of strata, corresponding to the subspaces in $\LL$ in the 
chosen linear extension order. 

Let us mention here that there is a projective
analogue~$\overline{Y}_{\AA}$ of the affine arrangement
model~$Y_{\AA}$ (cf.~\cite[\S 4]{DP1}). In fact, the affine model
$Y_{\AA}$ is the total space of a line bundle
over~$\overline{Y}_{\AA}$. A description of the integral cohomology
algebra of $\overline{Y}_{\AA}$ for complex hyperplane arrangements $\AA$ in
\cite{DP2} gave rise to a class of abstract algebras defined by atomic
lattices that bear a wealth of geometric meanings (cf.\ \cite{FY}).

We will need to refer to projective arrangement models only in one of our 
examples in Section~\ref{sec_appl}. We therefore stay with the affine setting in the 
following exposition.

\subsubsection{Normal crossing divisors and nested set stratification}
\label{sssec_stratification}
 The term {\em wonderful \/} models has been coined for $Y_{\AA}$ 
and its generalizations for other choices of building sets. We summarize 
the key facts about the maximal model supporting this connotation.

The space $Y_{\AA}$ is a smooth algebraic variety with a natural
projection onto the original ambient space $V$, $p:\, Y_{\AA} \lra V$.
The map $p$ is the projection onto the first coordinate of the ambient
space of $Y_{\AA}$ on the right hand side of~(\ref{eq_modelconstr}),
respectively the concatenation of blowdown maps of the sequence of
blowups resulting in~$Y_{\AA}$. This projection is an
isomorphism on $\MM(\AA)$, while the complement
$Y_{\AA}\,{\setminus}\,\MM(\AA)$ is a divisor with normal crossings
with irreducible components indexed by the elements of
$\LL_{>\hat 0}$. An intersection of several irreducible components is
non-empty (moreover, transversal and irreducible) if and only if the
indexing lattice elements form a totally ordered set, i.e., a chain, 
in~$\LL$~\cite[3.1,3.2]{DP1}. The stratification by irreducible components
of the divisor and their intersections is called the {\em nested set
  stratification\/} of $Y_{\AA}$, denoted $(Y_{\AA}, \mathfrak D)$,
for reasons that lie in the more general model construction for arbitrary
building sets rather than the maximal building set~$\LL_{>\hat 0}$.

\subsubsection{An encoding of points in maximal arrangement models}
\label{sssec_encoding} 
Points in $Y_{\AA}$ can be described as a sequence of a point and a number of
lines in the vector space~$V$ according to the form of the ambient
space for $Y_{\AA}$ given on the right hand side of~(\ref{eq_modelconstr}).
However, there is a lot of redundant information in that description. 
The following compact encoding of points was suggested in~\cite[Sect 4.1]{FK2}.

\begin{prop} \label{prop_encoding}
Let $\omega$ be a point in the maximal wonderful model $Y_{\AA}$ for 
a subspace arrangement $\AA$ in complex or real space~$V$.  Then $\omega$
can be uniquely written as
\begin{equation}\label{eq_encoding}
\omega \,\, = \,\, (x,H_1, \ell_1, H_2, \ell_2, \ldots, H_t,\ell_t)
        \,\, = \,\,(x,\ell_1, \ell_2, \ldots,\ell_t)\, , 
\end{equation}
where $x$ is a point in~$V$, the $H_1,\ldots,H_t$ form a descending chain 
of  subspaces in $\LL_{>\hat 0}$, and the~$\ell_i$ are lines in $V$, all 
subject to a number of additional conditions. 
\end{prop}

More specifically, $x\,{=}\,p(\omega)$, and the linear space $H_1$ is the 
maximal lattice element that, as a subspace of~$V$,
contains~$x$. The line $\ell_1$ is orthogonal to $H_1$ and
corresponds to the coordinate entry of $\omega$ indexed by $H_1$
in $\P(V/H_1)$.  The lattice element $H_2$, in turn, is the
maximal lattice element that contains both $H_1$ and $\ell_1$. The
specification of lines~$\ell_i$, i.e., lines that correspond to
coordinates of $\omega$ in $\P(V/H_i)$, and the construction of
lattice elements $H_{i+1}$, continues analogously for $i\geq 2$
until a last line $\ell_t$ is reached whose span with $H_t$ is not
contained in any lattice element other than the full ambient
space~$V$. Note that, if $H_t$ is a~hyperplane, then the line
$\ell_t$ is uniquely determined. The whole space $V$ can be
thought of as $H_{t+1}$.
Observe that the $H_i$ are determined by $x$
and the sequence of lines $\ell_i$; we choose to include the $H_i$ 
at times in order to keep the notation more transparent.

The full coordinate information on $\omega$ can be recovered from 
(\ref{eq_encoding}) by setting $H_0=\bigcap \AA$, $\ell_0=\langle x\rangle$, 
and retrieving the coordinate $\omega_{H}$ indexed by $H\,{\in}\,\LL_{> \hat 0}$ as 
 \begin{equation} \label{eq_recencoding}
  \omega_H\, \, = \, \, \langle \ell_j,H\rangle / H  \, \, \in \, \, \P(V/H)\, ,
\end{equation}
where $j$ is chosen from $\{1,\ldots, t\}$ such that
$H\leq H_j$, but $H\not \leq H_{j+1}$.
        
For completeness, let us mention here that we can tell the open stratum
in the nested set stratification $(Y_{\AA},\mathfrak D)$ that contains a 
given point $\omega$ from its point/line encoding stated in 
Proposition~\ref{prop_encoding}. 

\begin{prop}
{\rm (\cite[Prop 4.5]{FK2})}
A point $\omega$ in a maximal arrangement model $Y_{\AA}$ is contained 
in the open stratum of $(Y_{\AA},\mathfrak D)$ indexed by the chain
$H_1\,{>}\,H_2\,{>}\,\ldots \,{>}\,H_t\,{>}\,\hat 0$ in~$\LL$ if 
and only if its point/line description~{\rm (\ref{eq_encoding})} reads 
$\omega\,{=}\,(x,H_1, \ell_1, H_2, \ell_2, \ldots, H_t,\ell_t)$.

\end{prop}


\subsection{Group actions on arrangement models and a description 
            of stabilizers}\label{ssec_Garrgtmodels}

Provided an arrangement is invariant under the action of a 
finite group, this action extends to the maximal arrangement model. 
We review the details, and recall a description for stabilizers of
points in the model from~\cite{FK2}. 

\subsubsection{Group actions on $Y_{\AA}$} \label{sssec_Garrgtmodels}
Let $\AA$ be an arrangement that is invariant under the linear action 
of a finite group $G$ on the real or complex ambient space~$V$. Without
loss of generality, we can assume that this action is 
orthogonal~\cite[2.3,\,\,Thm~1]{V}. We denote the corresponding $G$-invariant
positive definite symmetric bilinear form by the usual scalar product. 

The group $G$ acts on the ambient space of the arrangement model $Y_{\AA}$,
i.e., for $(x,(x_{X})_{X\in \LL_{>\hat 0}})
      \,{\in}\,V\,{\times}\,\prod_{X\in\LL_{>\hat 0}}\, \P(V/X)$ and
$g\in G$, we have 
\begin{eqnarray*}
    g \,(x,(x_X)_{X\in \LL_{>\hat 0}}) & = &
   (\,g(x), (\,g(x_{g^{-1}(X)})\,)_{X\in\LL_{>\hat 0}})\, ,
\end{eqnarray*}
where $g(x_{g^{-1}(X)})\,{\in}\,\P(V/X)$ for $X\,{\in}\,\LL_{>\hat 0}$. 
Since the inclusion map $\Psi$ of~(\ref{eq_modelconstr}) commutes with the
$G$-action, and $G$ acts continuously on $V$, we conclude that
$Y_{\AA}=\overline{{\rm Im} \Psi}$ is as well  $G$-invariant. 
In particular, the  $G$-action on~$Y_{\AA}$ extends the $G$-action 
on the complement of $\AA$.

\subsubsection{Stabilizers of points on $Y_{\AA}$} \label{sssec_stabilizers}

The point/line description for points in the arrangement model $Y_{\AA}$
given in~\ref{sssec_encoding}  allows for a concise description of 
stabilizers with respect to the $G$-action on $Y_{\AA}$. 

\begin{prop}{\rm (\cite[Prop 4.2]{FK2})} \label{prop_stab}
  For a maximal arrangement model $Y_{\AA}$ that is  equipped with the action
  of a finite group $G$ stemming from a linear action of $G$ on the
  arrangement, the stabilizer of a point $\omega\,{=}\,(x,H_1, \ell_1, H_2,
  \ell_2, \ldots, H_t,\ell_t)$ in $Y_{\AA}$
  is of the form 
\begin{equation}\label{eq_stab}
\stab_{Y_{\AA}} (\omega) \, \, = \, \, \stab_V(x)\, \cap \,
                               \stab_V(\ell_1) \, \cap \, \ldots
                               \, \cap \, \stab_V(\ell_t)\, ,
\end{equation}
where, for  $i\,{=}\,1,\dots,t$, $\stab_V(\ell_i)$ denotes the elements in $G$ that
preserve the line $\ell_i$ in $V$ as a set.
\end{prop}


\subsection{Models for local subspace arrangements}\label{ssec_lamodels}

The arrangement model construction of De~Concini \& Procesi generalizes 
to the context of local subspace arrangements. 

\begin{df} \label{df_locsubsparrgts}
  Let $X$ be a smooth $d$-dimensional real or complex manifold and
  $\AA$ a family of smooth real or complex submanifolds in $X$ such
  that all non-empty intersections of submanifolds in $\AA$ are
  connected, smooth submanifolds. The family $\AA$ is called a {\em
    local subspace arrangement\/} if for any $x\,{\in}\, \bigcup \AA$
  there exists an open neighborhood $U$ of $x$ in~$X$, a subspace
  arrangement $\tilde \AA$ in a real or complex $d$-dimensional vector
  space $V$ and a diffeomorphism $\phi:\, U\lra V$, mapping $\AA$ to
  $\tilde \AA$.
\end{df}

Local subspace arrangements fall into the class of conically
stratified manifolds as appearing in work of MacPherson \&
Procesi~\cite{MP} in the complex and in work of Gaiffi~\cite{Ga} 
in the real setting.

A generalization of the arrangement model construction of De~Concini \&
Procesi by sequences of blowups of smooth strata for conically
stratified complex manifolds is given in~\cite{MP}. Details are
provided for blowing up so-called irreducible strata, the more general
construction for an arbitrary building set of the stratification is
outlined in Sect.~4 of~\cite{MP}. 

In this article, we will be concerned with {\em maximal 
wonderful models\/}  for conically stratified real manifolds~$X$, 
in the special case of local subspace arrangements~$\AA$. 
The maximal model $Y_{\AA}\,{=}\,Y_{X,\AA}$ results
from successive blowups of {\em all\/} initial strata, respectively 
their proper transforms, according to some linear order on strata which 
is non-decreasing in dimension.

In fact, local subspace arrangements consisting of a {\em finite\/}
number of submanifolds implicitly appear already in the arrangement
model construction of De~Concini \& Procesi~\cite{DP1}.  A single
blowup in a subspace arrangement leads to the class of local
arrangements, and it is due to the choice of blowup order on building
set strata that this class is closed under blowups that occur in the
inductive construction of the arrangement models (cf.\ the
discussion in~\cite[4.1.2]{FK1}, in particular, Example 4.6).

We will encounter the case of local subspace arrangements $\AA$ in 
a smooth real manifold~$X$ 
that are invariant under the differentiable action of a finite group~$G$ 
on $X$. The $G$-action can be extended to the maximal model $Y_{\AA}$, 
observing that we can simultaneously blow up orbits of strata, thereby
lifting the $G$-action step by step through the construction process. 
In particular, the concatenation of blowdown maps $p:\,Y_{\AA}\ra X$ 
is $G$-equivariant.


\section{Digitalizing finite group actions}
\label{sec_digital}


\subsection{Finite linear actions on $\R^n$}
\label{ssec_lin}

\noindent
In this subsection we assume $G$ to be a finite subgroup of the orthogonal 
group $O(n)$ acting
effectively on $\R^n$. As pointed out before, assuming the action to be 
orthogonal is not a restriction (cf.\ \ref{sssec_Garrgtmodels}).

We construct an abelianization of the given action. For any subgroup
$H$ in $G$ (we use the notation $H\,{\leq}\,G$ in the sequel), define 
\[
   L(H)\, \, := \, \, \langle \, \ell\, |\, 
                      \ell \mbox{ line in }\, \R^n \mbox{ with }\,
                      h\cdot \ell=\ell \mbox{ for all }\, h\in H\,\rangle \, ,
\]
the linear span of lines in $\R^n$ that are invariant under~$H$, i.e.,
the span of lines that are either fixed or flipped by any element
$h$ in~$H$. Denote by $\AA$ the arrangement given by the {\em proper\/}
subspaces $L(H)\,{\subsetneq}\,\R^n$, $H$ subgroup in $G$. 
Set $Y_{\AA}$ denote the maximal De~Concini-Procesi wonderful model for $\AA$
as discussed in~\ref{ssec_arrgtmodels}. If we want to stress the particular 
group action that gives rise to the arrangement~$\AA$ we write $\AA(G)$
and $Y_{\AA(G)}$, or $\AA(G\car\R^n)$ 
and $Y_{\AA(G\car\R^n)}$,
respectively. 

We will now propose $Y_{\AA}$ as an abelianization of the given linear action.
Recall that we use the term {\em digitalization\/} for an abelianization
with stabilizers that are not merely abelian but elementary
abelian $2$-groups, i.e., are isomorphic to $\Z_2^k$ for some $k\,{\in}\,\N$.

\begin{thm} \label{thm_dglin}
Let an effective action of a finite subgroup $G$ of $O(n)$ on $\R^n$ be given.
Then the wonderful arrangement model $Y_{\AA(G)}$, as defined above, is a 
digitalization of the given action.  
\end{thm}

\begin{pf}
As a first step we prove that
\[
                 L({\rm stab}\,\omega) \, \, =\, \, \R^n\, ,\qquad 
                                    \mbox{ for any }\, \omega\in Y_{\AA}\, . 
\]
Let $\omega\in Y_{\AA}$. Using the encoding of points in arrangement models as
sequences of point and lines from~\ref{sssec_encoding}, we have 
$\omega\,{=}\,(x,\ell_1,\ldots,\ell_t)$, the associated sequence of building 
set spaces being $V_1,\ldots,V_t$. The description of stab$\,\omega$ 
from Proposition~\ref{prop_stab},
\[
   {\rm stab}\,\omega \, \, = \, \, 
   {\rm stab}\,x \, \cap \, {\rm stab}\,\ell_1 \, \cap \, \ldots 
                                        \, \cap \, {\rm stab}\,\ell_t\, , 
\]
implies that $x\,{\in}\,L({\rm stab}\,\omega)$, and 
$\ell_i\,{\subseteq}\,L({\rm stab}\,\omega)$ for $i\,{=}\,1,\ldots,t$.

The building set element $V_1$ is the smallest subspace among intersections 
of spaces $L(H)$ in $\AA$ such that $x\,{\in}\,V_1$, in particular,
$V_1\,{\subseteq}\,L({\rm stab}\,\omega)$. Similarly, the building set element 
$V_2$ is the smallest subspace among intersections 
of spaces $L(H)$ in $\AA$ such that $\langle V_1,\ell_1\rangle\,{\subseteq}\,V_2$;
since $\langle V_1,\ell_1\rangle\,{\subseteq}\,L({\rm stab}\,\omega)$, 
so is~$V_2$: $V_2\,{\subseteq}\,L({\rm stab}\,\omega)$.  

By analogous arguments we conclude that $V_3,\ldots,V_{t+1}\,{\subseteq}\,L({\rm
  stab}\,\omega)$. However, by the description of
$\omega$ as a sequence of point and lines we know that
$V_{t+1}\,{=}\,\R^n$, which proves our claim.

\smallskip
With  $L({\rm stab}\,\omega)\,{=}\,\R^n$, we can now choose a basis 
$v_1,\ldots, v_n$ in $\R^n$ such that any $\langle v_i\rangle$, for $i=1,\ldots,n$,
is invariant under the action of $\mathrm{stab}\, \omega$.

Consider the homomorphism
\begin{eqnarray*}
  \alpha \,: \quad {\rm stab}\,\omega & \longrightarrow & \Z_2^n \\
                                h  & \longmapsto     & 
                                     (\epsilon_1,\ldots,\epsilon_n)\, , 
\end{eqnarray*}
with $\epsilon_i\,{\in}\,\Z_2$ defined by $h(v_i)\,{=}\,\epsilon_i\, v_i$ 
for $i\,{=}\,1,\ldots,n$.
Since we assume the action to be effective, $\alpha$ is injective. Hence
stab$\,\omega\,{\cong}\,\Z_2^k$ for some $k\,{\leq}\,n$. 
\end{pf}

 \subsection{Finite differentiable actions on manifolds}
\label{ssec_mnfds}

\noindent
We now generalize the results of the previous subsection to differentiable
actions of finite groups on smooth manifolds. To this end, we first
propose a stratification of the manifold and show that the
stratification locally coincides with the arrangement stratifications
on tangent spaces that arise from the induced linear actions as described
in the previous section. We can assume, without loss of generality, that
the manifold is connected, since we can work with connected components 
one at a time.

\subsubsection{The $\LL$-stratification}
Let $X$ be a smooth manifold, $G$ a finite group that acts
differentiably on~$X$. For any point $x\,{\in}\,X$, and any subgroup
$H\,{\leq}\,{\rm stab}\,x$, $H$ acts linearly on the tangent space $T_xX$ of $X$ in
$x$. Consider as above
\[
   L(x,H)\, \, := \, \, \langle \, \ell\, |\, 
                      \ell \mbox{ line in }\, T_xX \mbox{ with }\,
                      h\cdot \ell=\ell \mbox{ for all }\, h\in H\,\rangle \, ,
\]
the linear subspace in $T_xX$ spanned by lines that are invariant under the 
action of~$H$. Denote the arrangement of proper subspaces $L(x,H)$ in $T_xX$, 
$\AA({\rm stab}\, x \car T_xX)$, by $\AA_x$.

For any subgroup $H$ in ${\rm stab}\, x$, we take up the homomorphism that 
occurred in the proof of Theorem~\ref{thm_dglin}, and define
\[
     \alpha_{x,H}:\,\, H \,\, \longrightarrow \,\, \Z_2^{\mathrm{dim}\,L(x,H)}
\]
by choosing a basis $v_1,\ldots,v_t$, $t\,{:=}\,\mathrm{dim}\,L(x,H)$, for $L(x,H)$, 
and setting
\[
    \alpha_{x,H} (h) \,=\, (\epsilon_1,\ldots,\epsilon_t)\, ,
\]
for $h\,{\in}\,H$, with $\epsilon_i\,{\in}\,\Z_2$ determined by $h(v_i)\,{=}\,\epsilon_i v_i$
for $i\,{=}\,1,\ldots,t$.

Moreover, we define 
\[
      F(x,H)\, := \, \mathrm{ker}\,  \alpha_{x,H}\, . 
\] 
Note that  $F(x,H)$ is the normal subgroup of elements in $H$ that fix 
all of $L(x,H)$ point-wise.
We denote by $\LL(x,H)$ the connected component of $\mathrm{Fix}(F(x,H)\,{\car}\,X)$ in $X$
that contains~$x$. 

Consider the stratification of $X$ by the collection
of submanifolds $\LL(x,H)$ for $x\,{\in}\,X$, $H\,{\leq}\,\mathrm{stab}\,x$,
\[
        \LL \, \, := \, \,\{\LL(x,H)\}_{x\in X, \, H \leq {\rm stab}\,x }\, .
\]
We will refer to this stratification as the {\em $\LL$-stratification\/}
of~$X$. Observe that $\LL$ is a locally finite stratification. 

\smallskip 

We recall the following fact from the theory of group actions on smooth 
manifolds:

\begin{prop} \label{prop_slice}
Let $G$ be a compact Lie group acting differentiably on a smooth 
manifold~$X$, and let $x_0\,{\in}\,X$. Then there exists a 
{\rm stab}$\,x_0$-equivariant diffeomorphism $\Phi_{x_0}$ from an open 
neighborhood
 $U$ of $x_0$ in $X$ to the tangent space $T_{x_0}X$ of $X$ in $x_0$.
\end{prop}
 
This is a special case of the so-called {\em slice theorem\/}~\cite{A,tD}
that originally appeared in work of Bochner~\cite{Bo}.

We return to our setting of $G$ being a finite group.
 
\begin{prop} \label{prop_mapstrat}
The diffeomorphism $\Phi_{x_0}$ maps the $\LL$-stratification of $X$ to the
arrangement stratification on $T_{x_0}X$ given by $\AA_{x_0}$, i.e.,
\[
   \Phi_{x_0}(\LL(x_0,H))\, \, = \, \, L(x_0,H) \qquad 
                           \mbox{ for any }\, H\leq {\rm stab}\,x_0\,. 
\]
\end{prop}

\begin{pf}
  By definition, $\LL(x_0,H)\,{=}\,{\rm Fix}\,(F(x_0,H)\car X)$,
  which, using the stab$\,x_0$-equivariance of $\Phi_{x_0}$, implies
  that $\Phi_{x_0}(\LL(x_0,H))\,{=}\, {\rm Fix}\,(F(x_0,H) \car
  T_{x_0}X)$. We are left to show that
\[
   {\rm Fix}\,(F(x_0,H) \car T_{x_0}X) \, = \, L(x_0,H)\, .                
\]

Obviously, $L(x_0,H)\,{\subseteq}\,{\rm Fix}\,(F(x_0,H)\car
T_{x_0}X)$, and we need to see that Fix$\,(F(x_0,H)\car T_{x_0}X)$
does not exceed~$L(x_0,H)$.

Note that $H$ acts on $L(x_0,H)$. By definition, $F(x_0,H)$ is a
normal subgroup of $H$ with quotient $H/F(x_0,H)\,{\cong}\,\Z_2^d$ for
some $d\,{\leq}\,t\,{=}\,\mathrm{dim}\, L(x_0,H)$, and we find that
$H$ acts on Fix$\,(F(x_0,H) \car T_{x_0}X)$: For $x\,{\in}\,{\rm
  Fix}\,(F(x_0,H) \car T_{x_0}X)$, $h\,{\in}\,H$, and
$h_1\,{\in}\,F(x_0,H)$, we have $h_1hx\,{=}\,h\tilde h_1x$ for some
$\tilde h_1\,{\in}\,F(x_0,H)$, thus $h_1hx\,{=}\,hx$, i.e.,
$hx\,{\in}\,{\rm Fix}\,(F(x_0,H) \car T_{x_0}X)$.

Instead of considering the action of $H$ on ${\rm Fix}\,(F(x_0,H) \car
T_{x_0}X)$, we consider the induced action of $H/F(x_0,H)$ on ${\rm
  Fix}\,(F(x_0,H) \car T_{x_0}X)$.  Since
$H/F(x_0,H)\,{\cong}\,\Z_2^d$ for some $d\,{\leq}\,t$, ${\rm
  Fix}\,(F(x_0,H) \car T_{x_0}X)$ decomposes into $1$-dimensional
representation spaces, which, as lines that are invariant under the
action of $H$, must be contained in $L(x_0,H)$ by definition. This
shows that ${\rm Fix}\,(F(x_0,H) \car T_{x_0}X)$ does not exceed
$L(x_0,H)$, and thus completes our proof.
\end{pf}

\begin{remrm} \label{rem_linear}
Applying Proposition~\ref{prop_mapstrat} to a linear action $G\,{\car}\,X\,{=}\,\R^n$
for $x_0\,{=}\,0$, we see that the linear subspaces $L(H)$, for $H\,{\leq}\,G$,
are fixed point sets of subgroups of~$H$, namely
$L(H)\,{=}\,\fix(F(H))$,
where $F(H)$ is the subgroup of elements in $H$ that fix $L(H)$
point-wise, $F(H)\,{=}\,\{h\,{\in}\, H\,|\,hx=x \mbox{ for all }\,
x\in L(H)\}$.
\end{remrm}

In particular, Proposition~\ref{prop_mapstrat} shows that the
submanifolds $\LL(x,H)$ in the $\LL$-stratification form a local
subspace arrangement in $X$.  Moreover, the $\LL$-stratification is
invariant under the action of $G$ since $g(\LL(x,H))\,{=}\,\LL(g(x),
gHg^{-1})$ for any $x\,{\in}\,X$, $H\,{\leq}\,{\rm stab}\,x$, and any
$g\,{\in}\,G$. Hence, we have at hand the maximal $G$-equivariant 
wonderful model $Y_{\LL}\,{=}\,Y_{X,\LL}$ of the local subspace 
arrangement $\LL$ in $X$ as outlined in~\ref{ssec_lamodels}.

\subsubsection{Digitalizing manifolds}

We propose the maximal wonderful model of $X$ with respect to
the \mbox{$\LL$-stra}tification as a digitalization of the manifold~$X$.

\begin{thm} \label{thm_dgman}
  Let $G$ be a finite group acting differentiably and effectively
  on a smooth manifold~$X$. Then the maximal wonderful model of $X$ with
  respect to the \mbox{$\LL$-stra}tification $Y_{X,\LL}$ is
  a digitalization of the given action.
\end{thm}

\begin{pf}
Let $x$ be a point in $Y_{X,\LL}$, $x_0\,{=}\,p(x)$ its image under the 
blowdown map $p:\, Y_{X,\LL} \lra X$. Since $p$ is $G$-equivariant,
stab$\,x\,{\subseteq}\,{\rm stab}\,x_0$, hence we can restrict our 
attention to stab$\,x_0$ when determining the stabilizer of $x$ in~$G$.

Consider the stab$\,x_0$-equivariant diffeomorphism $\Phi_{x_0}$
as discussed above (Proposition~\ref{prop_slice}),
\[
\Phi_{x_0}\,: \quad U \, \, \lra \, \, T_{x_0}X\, ,
\]
where
$U$ is an open neighborhood of $x_0$ in $X$, such that $\Phi_{x_0}$ maps 
the $\LL$-stratification on~$U$ to the arrangement stratification on the
tangent space at $x_0$. Since the De~Concini-Procesi model is defined
locally, the diffeomorphism $\Phi_{x_0}$ induces a
stab$\,x_0$-equivariant diffeomorphism between the inverse image of
$U$ under the blowdown map, $p^{-1}U\,{=}\,Y_{U,\LL}$, and the
De~Concini-Procesi model for the arrangement $\AA_{x_0}$ in the tangent
space $T_{x_0}X$,
\[
\widetilde \Phi_{x_0}\,: 
 \quad Y_{U,\LL} \, \, \lra \, \, Y_{T_{x_0}X, \AA_{x_0}}\, .
\]
In particular, 
\[
    {\rm stab}\,x \, \, \cong \, \, {\rm stab}\, \widetilde \Phi_{x_0}(x)\, ,
\]
which, by our analysis of the linear setting, is an elementary abelian $2$-group, 
provided we can see that stab$\,x_0$ acts
effectively on $T_{x_0}X$. To settle this remaining point, assume that
there exists a group element $g\,{\neq}\,e$ in stab$\,x_0$ that fixes
all of $T_{x_0}X$. By Proposition~\ref{prop_slice}, $g$ then fixes an
open neighborhood of $x_0$ in~$X$, which implies that $g$ fixes all of
$X$, contrary to our assumption of the action being effective.
\end{pf}


\section{Permutation actions on linear and on projective spaces}
\label{sec_appl}

One of the most natural linear actions of a finite group is the action
of the symmetric group $\Ss_n$ permuting the coordinates of a real
$n$-dimensional vector space. This action induces a differentiable
action of $\Ss_n$ on $(n{-}1)-$dimensional real projective
space~$\R\P^{n-1}$. Our goal in this section is to give explicit
descriptions of the $\LL$-stratifications and the resulting
digitalizations in both cases. The {\em braid arrangement\/} $\AA_{n-1}$
will play a central role for both stratifications; we thus recall that 
$\AA_{n-1}\,{=}\, \{H_{ij}{:}\, x_j{-}x_i{=}0\,|\,1{\leq} i,j{\leq}n\}$.
Its intersection lattice is the lattice $\Pi_n$ of set partitions 
of~$[n]:=\{1,\ldots,n\}$ ordered by reversed refinement.

We will show that, in the case of the real permutation action, the
arrangement $\AA(\Ss_n)$ coincides with the rank~$2$ truncation of the
braid arrangement, $\AA_{n-1}^{\mathrm{rk}\geq 2}$, i.e., the braid
arrangement $\AA_{n-1}$ without its hyperplanes.  The abelianization
construction proposed in the present article hence specializes to the
maximal model of the braid arrangement discussed in \cite{FK2}.

For the permutation action on~$\R\P^{n-1}$, we show that the
$\LL$-stratification coincides with the rank $2$~truncation of the
projectivized braid arrangement, $\P\AA_{n-1}^{\mathrm{rk}\geq 2}$,
thus the digitalization proposed in~\ref{ssec_mnfds} coincides with
the maximal projective arrangement model for~$\AA_{n-1}$ (cf.~\cite[\S
4]{DP1}).

Let us fix some notation that will come in handy when dealing with the
permutation action on $\R^n$, and the induced action on $\R\P^{n-1}$,
respectively. Any permutation $\sigma\in \Ss_n$ determines a set
partition of $[n]$, $\rho(\sigma)\,{\vdash}\,[n]$, by
its cycle decomposition. Moreover, any set partition
$\pi\,{=}\,(B_1|\ldots|B_k)\,{\vdash}\,[n]$ gives rise to an
intersection of hyperplanes $U_{\pi}$ in the braid arrangement
$\AA_{n-1}$, namely
\[
  U_{\pi} \, \, := \, \, \bigcap_{r=1}^k \, \bigcap_{i,j\in B_r}\, 
                          H_{ij}\, .    
\] 
We call  $U_{\pi}$ the {\em braid space} associated to~$\pi$.


\subsection{Digitalizing the real permutation action}
\label{ssec_MaxAn}

As outlined above, we will now recover the rank $2$~truncation of the braid 
arrangement as the arrangement $\AA(\Ss_n)$ arising from the real permutation 
action.     

\begin{thm} \label{thm_ASn}
  The arrangement $\AA(\Ss_n)$ associated with the real permutation
  action as described in~{\rm \ref{ssec_lin}} coincides with the rank
  $2$~truncation of the braid arrangement. In particular, the
  digitalization $Y_{\AA(\Ss_n)}$ of {\rm Theorem~\ref{thm_dglin}}
  coincides with the maximal wonderful model of the braid arrangement
  as discussed in \cite{FK2}.
\end{thm}

\begin{pf}
  We first show that any proper subspace $L(H)$ in $\R^n$, for $H$ a
  subgroup of $\Ss_n$, is a braid space of codimension at
  least~$2$:\newline Recall from Remark~\ref{rem_linear} that any
  subspace $L(H)$ is a fixed point set, namely,
  $L(H)\,{=}\,{\fix}(F(H))$ where $F(H)$ is the subgroup of elements
  in $H$ that fix $L(H)$ point-wise.  Obviously,
  $\fix(F(H))=\cap_{\sigma \in F(H)}\fix(\sigma)$. Therefore, it is
  enough to prove that $\fix(\sigma)$ is a braid space for any
  permutation $\sigma\in\Ss_n$. However, it is obvious that
  $\fix(\sigma){=}U_{\rho(\sigma)}$, where $\rho(\sigma)$ is the cycle
  decomposition of $\sigma$ as defined above.  Observe here that proper
  subspaces of type $L(H)$ are never of codimension~$1$; if $L(H)$
  were that, then $H$ would also leave the orthogonal line $L(H)^{\perp}$
  invariant, in contradiction to $L(H)$ being proper.
  
  We will now prove that all braid spaces other than hyperplanes occur
  in the arrangement $\AA(\Ss_n)$. To this end, we show that we can
  realize braid spaces $U_{\pi}$, $\pi\,{\vdash}\,n$, with
  type$(\pi){=}(3,1^{n-3})$ and type$(\pi){=}(2^2,1^{n-4})$ as
  subspaces $L(H)$ for some subgroups $H$ of~$\Ss_n$. Any braid space
  of higher codimension will be obtained as an intersection of those.

  Without loss of generality, set $\pi\,{=}\,123|4|\ldots|n$ to cover 
  elements of the first type. 
  Let~$H$ be the Young subgroup $\Ss_{\{1,2,3\}}\times\Ss_{\{4\}}\times 
  \dots \times\Ss_{\{n\}}$. Clearly, if a line is invariant under the 
  action of $H$, then it must be fixed point-wise. This is the case if 
  and only if the first~$3$ coordinates of a generating vector are equal. 
  We conclude that $L(\Ss_{\{1,2,3\}})\,{=}\,U_{123}$.
  
  Now, set $\pi\,{=}\,12|34|5|\ldots|n$.  Let $H$ be the subgroup of
  $\Ss_n$ generated by the transpositions $(12)$ and $(34)$, and the
  involution $(13)(24)$. We see that $H$ is isomorphic to a~wreath
  product $\Ss_2\wr\Ss_2$. The subspace of points fixed by $H$ is the
  braid space~$U_{1234}$. Furthermore, the line spanned by the
  vector $(1,1,-1,-1,0,\dots,0)$ belongs to $L(H)$, it is flipped
  by~$H$.  It is not difficult to see that~$U_{1234}$ and the line
  $\langle (1,1,-1,-1,0,\dots,0)\rangle$ span the entire~$L(H)$; on the 
  other hand, they span~$U_{12|34}$, so we
  conclude that $L(\Ss_2\wr\Ss_2)\,{=}\,U_{12|34}$. \newline
\mbox{ }
\end{pf}


\subsection{Digitalizing the permutation action on real projective space}
\label{ssec_MaxPAn} 

We will consider the $\LL$-stratification on $\R\P^{n-1}$ induced by the 
permutation action of $\Ss_n$, and we will give a description of the 
digitalization proposed in~\ref{ssec_mnfds}.

Let us first prove a lemma that describes the fixed point set of a
single permutation $\sigma\,{\in}\,\Ss_n$ on $\R\P^{n-1}$. The fact
that $\LL$-strata are defined as connected components of fixed point
sets, points to the significance of this lemma for the following
considerations.

\begin{lm} \label{lm_fixsigma}
Let $\sigma$ be a permutation in $\Ss_n$. The fixed point set of $\sigma$ on
$\R\P^{n-1}$ decomposes as a disjoint union of submanifolds
\begin{equation} \label{eq_fixsigma}
    \fix(\sigma \car \R\P^{n-1})\,\, = \,\, 
                      \P U_{\rho(\sigma)} \, \cup \, 
                      \P V_{\sigma}\, , 
\end{equation}
where $V_{\sigma}$ is the span of lines generated by vectors $v_B$,
one for each even length block $B$ in the partition~$\rho(\sigma)$
associated with~$\sigma$. The vector $v_B$ has sign-alternating
entries $\pm 1$ in the coordinates contained in~$B$ in the linear
order determined by the iteration of~$\sigma$, and has entries~$0$
otherwise.
\end{lm}

\begin{pf}
  Let $L^+(\sigma)$ denote the span of lines in $\R^n$ that are
  point-wise fixed by~$\sigma\,{\in}\,\Ss_n$ with respect to the
  permutation action on $\R^n$, and let $L^-(\sigma)$ denote the span
  of lines in~$\R^n$ that are flipped by~$\sigma$. Obviously,
\[
  \fix(\sigma \car \R\P^{n-1})\,\, = \,\, \P(L^+(\sigma)\,\cup\, L^-(\sigma))\, .   
\]
First we note that $L^+(\sigma)\,=\, U_{\rho(\sigma)}$. To describe
$L^-(\sigma)$ we proceed as follows: In order for a line $l=\langle
v\rangle$, $v\in\R^n$, to be flipped by the action of $\sigma$, the
non-zero coordinates of~$v$ must have the same absolute value within
each cycle of $\sigma$, i.e., each block of $\rho(\sigma)$, with
actual entries alternating in sign in the order prescribed by
iteration of~$\sigma$. This implies that the coordinates
of~$v$, which are contained in the odd length blocks of $\rho(\sigma)$,
must be~$0$. Moreover, we see that $L^-(\sigma)$ is
generated by vectors $v_B$ for the even length blocks $B$ in
$\rho(\sigma)$, with $v_B$ having alternating entries $\pm 1$ in $B$
and $0$ otherwise. Thus, $L^-(\sigma)\,{=}\,V_{\sigma}$ as described in the 
statement of the Lemma. The subspaces
$L^+(\sigma)$ and~$L^-(\sigma)$ are orthogonal in $\R^n$, hence
\[
    \fix(\sigma \car \R\P^{n-1})\, = 
    \,\P(U_{\rho(\sigma)}\,\cup\, V_{\sigma})\, = 
    \,\P U_{\rho(\sigma)}\,\cup\, \P V_{\sigma} \, , 
\]
and the union is disjoint. 
\end{pf}

\begin{thm} \label{thm_LPn}
The $\LL$-stratification on $\R\P^{n-1}$ induced by the permutation action 
of $\Ss_n$ coincides with the rank $2$~truncation of the projectivized braid 
arrangement. In particular, the digitalization $Y_{\R\P^{n-1}, \LL}$
coincides with the maximal projective arrangement model for~$\P \AA_{n-1}$.  
\end{thm}

\begin{pf}
  We first show that any projectivized braid space of codimension at
  least~$2$ occurs as an $\LL$-stratum on $\R\P^{n-1}$. As in the proof of 
  Theorem~\ref{thm_ASn}, it will be
  enough to show that we can realize projectivizations of braid
  spaces $\P U_{\pi}$ with type$(\pi){=}(3,1^{n-3})$ and
  type$(\pi){=}(2^2,1^{n-4})$ as $\LL$-strata $\LL(\ell,H)$ for some points
  $\ell \,{\in}\,\R\P^{n-1}$ and subgroups $H$ of~$\Ss_n$.
  
  Again, we start by setting $\pi\,{=}\,123|4|\ldots |n$ without loss
  of generality.  We choose $\ell\,{=}\,\langle(1,\ldots,1)\rangle$ as
  the reference point in $\R\P^{n-1}$, and $H\,{=}\,\langle
  (123)\rangle\,{\cong}\, \Z_3$ as subgroup of~$\Ss_n$.  Observe that,
  due to our special choice of~$\ell$, $H$ acts on the tangent space
  $T_{\ell}\R\P^{n-1}$ by permuting coordinates.  We find that
  $L(\ell,\Z_3)\,{=}\,U_{123}$, and $F(\ell,\Z_3)\,{=}\,\Z_3$.
  Referring to Lemma~\ref{lm_fixsigma}, we see that $\fix(\Z_3 \car
  \R\P^{n-1})\,{=}\,\fix((123) \car \R\P^{n-1})$ has only one
  connected component $\P U_{123}$.  The second submanifold
  in~(\ref{eq_fixsigma}) does not come into play, since there are no
  even length block sizes occurring in the partition $123|4|\ldots
  |n$. We conclude that $\LL(\ell,H)\,=\,\P U_{123}$.
  
  We now set $\pi\,{=}\,12|34|5|\ldots |n$, and we choose
  $\ell\,{=}\,\langle(1,\ldots,1)\rangle$ in $\R\P^{n-1}$ as before, and $H\,{=}\,\langle
  (1234)\rangle\,{\cong}\, \Z_4$ as subgroup of $\Ss_n$. 
  Again, $H$ acts on $T_{\ell}\R\P^{n-1}$
  by permuting coordinates and we see that
  $L(\ell,\Z_3)\,{=}\,U_{12|34}$. The latter space is {\em spanned\/}
  by the lines that are fixed by $(1234)$, i.e., by
  $U_{1234}$, {\em and\/} the line $\langle(1,-1,1,-1,0,\ldots,0)\rangle$
  which is flipped by~$(1234)$. We find that $F(\ell,\Z_4)\,{=}\,
  \langle (13)(24)\rangle\,{\cong}\,\Z_2$, and 
  $\LL(\ell,\Z_4)$ is the connected component of 
  $\fix((12)(34)
  \car \R\P^{n-1})\,=\,\P U_{12|34} \cup \P V_{(12)(34)}$ that contains~$\ell$.
  The first submanifold contains our chosen point~$\ell$ and is disjoint from 
  the second, hence $\LL(\ell,\Z_4)\,{=}\,\P(U_{12|34})$. 
  
  Summarizing up to this point, we see that the rank~$2$ truncation of
  the projectivized braid arrangement is part of the
  $\LL$-stratification induced by the permutation action on real
  projective space. It remains to show that all $\LL$-strata indeed
  are projectivized braid spaces.
  
  To this end, let $\ell$ be a point in $\R\P^{n-1}$, i.e., 
  $\ell\,{=}\,\langle v \rangle$ is a line in 
  $\R^n$ with generating vector $v$ of unit length, and let $H$ be 
  a subgroup of stab$\,\ell$.

  We interpret the tangent space $T_{\ell}\R\P^{n{-}1}$ as the linear
  hyperplane orthogonal to $\langle v\rangle$ in $\R^n$,
  $T_{\ell}\R\P^{n{-}1}\,{=}\,\langle v\rangle^{\perp}\,{:=}\,T$.
  Instead of the linear stab$\,\ell$-action on the tangent space~$T$
  induced by $\Ss_n$ acting on $\R\P^{n-1}$, we can consider the
  permutation action of stab$\,\ell$ on~$T$. These two actions differ by at
  most a sign, which in particular implies that the construction of
  $L(H)$ yields the same subspaces with respect to both actions.

  We moreover observe that we can consider stab$\,\ell$ acting by permutation 
  of coordinates on $T\,{\oplus}\,\langle v \rangle \,{\cong}\,\R^n$ and we 
  obtain the space
  $L(H\,{\car}\,T )$ by restricting $L(H\,{\car}\,\R^n)$ to $T$. The latter
  space {\em is\/} a braid space of codimension at least $2$ by Theorem~\ref{thm_ASn}, 
\[
     L(H\,{\car}\,T ) \,\oplus \, \langle v \rangle \,\, = \, \, U_{\pi}    \, , 
\] 
  for some partition $\pi\,{\vdash}\,n $, rk$\pi\,{\geq}\,2$.

  Translating $T$ along $\ell$ into $T\,{+}\,v$, we can give an
  explicit description of the Bochner map $\Phi_{\ell}$ (cf.\ 
  Proposition~\ref{prop_slice}) that maps a neighborhood $U$ of $\ell$ in
  $\R\P^{n{-}1}$ differentiably and stab$\,\ell$-equivariantly to
  the tangent space $T_{\ell}\R\P^{n{-}1}$,
  \begin{eqnarray*} 
  \Phi_{\ell}\,: \qquad U & \lra & T_{\ell}\R\P^{n{-}1} \\
                     u  & \longmapsto & u \,\cap \,(\ell^{\perp}\,{+}\,v)\, . 
  \end{eqnarray*}
  We see that $L(H)\,{+}\,v$ in $T_{\ell}\R\P^{n{-}1}$ has the
  projectivization $\P U_{\pi}$ as its inverse image under~$\Phi_{\ell}$. 
  By Proposition~\ref{prop_mapstrat} we conclude that $\LL(\ell,H)$
  indeed is the projectivization of a braid space, which completes our
  proof.
\end{pf}

\begin{explrm} \label{ex_S3P}
  To illustrate our theorem on the $\LL$-stratification of real
  projective space induced by the permutation action and the resulting
  digitalization, we look at $\Ss_3$ acting on $\R\P^2$ in some
  detail.
  
  We depict $\R\P^2$ using the upper hemisphere model in
  Figure~\ref{fig_RP2}, where we place $\P\Delta^{\perp}$, for
  $\Delta\,{=}\,\langle (1,1,1)\rangle$, on the equator.

\begin{figure}[ht]
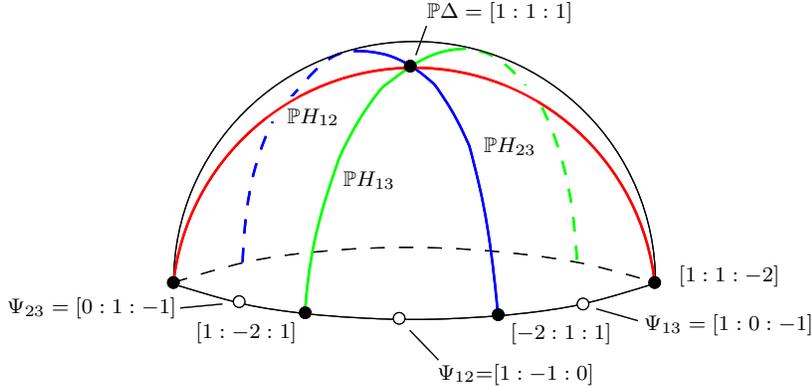

  \begin{picture}(0,0)%
    \includegraphics{RP2.pstex}%
  \end{picture}%
  \input{RP2.pstex_t}%
  
\caption{$\R\P^2$ stratified by loci of non-trivial stabilizers.}
\label{fig_RP2}
\end{figure}

The locus of points in $\R\P^2$ with non-trivial stabilizer groups
consists of the three lines~$\P H_{ij}$, $1\,{\leq}\,i,j\,{\leq}\,3$,
which are projectivizations of the hyperplanes in $\AA_2$,
intersecting in $\P \Delta\,{=}\,[1{:}1{:}1]$, and points $\Psi_{ij}$
on $\P\Delta^{\perp}$, where $\Psi_{ij}$ is the line orthogonal to
$H_{ij}$ in $\R^3$ for $1\,{\leq}\, i,j \,{\leq}\,3$.

Observe that the transposition $(i,j)\,{\in}\,\Ss_3$ acts on $\R\P^2$
as a central symmetry in $\Psi_{i,j}$, respectively, as a reflection
in $\P H_{i,j}$, which we illustrate in Figure~\ref{fig_action}.

\begin{figure}[ht]
  \begin{picture}(0,0)%
    \includegraphics{action.pstex}%
  \end{picture}%
  \input{action.pstex_t}%
  
\caption{$(i,j)\,{\in}\,\Ss_3$ acting on $\R\P^2$.}
\label{fig_action}
\end{figure}

We find that the arrangements $\AA_{\ell}({\rm
  stab}\,\ell\,{\car}\,T_{\ell}\R\P^2)$ associated with the induced
linear actions of the stabilizers on tangent spaces for $\ell\,{\in
}\,\R\P^2$ are empty unless $\ell\,{=}\,[1{:}1{:}1]$. In this case, we
see that $\Ss_3\,{\car}\, T_{[1{:}1{:}1]}\R\P^2$ coincides with the
standard action of $\Ss_3$ on $\R^3/ \Delta$, since transpositions, as
we observed above, act as reflections in the hyperplanes of the
projectivized braid arrangement. Thus,
$\AA_{[1{:}1{:}1]}(\Ss_3\,{\car}\,T_{[1{:}1{:}1]}\R\P^2)$ coincides
with the rank~$2$ truncation of the braid arrangement consisting of
the origin in the tangent space.

We conclude that the $\LL$-stratification is given by the single point
$[1{:}1{:}1]$ in $\R\P^2$, hence the digitalization we propose is the
blowup of $\R\P^2$ in this point,
\[
    Y_{\R\P^2,\LL}\, \, = \, \, {\rm Bl}_{[1{:}1{:}1]} (\R\P^2)\, .
\]
Topologically, this means to glue a M\"obius band into a pointed
$\R\P^2$, equivalently, to glue two M\"obius bands along their
boundaries, resulting in a Klein bottle (cf.\ Figure~\ref{fig_blowup}). 

\begin{figure}[ht]
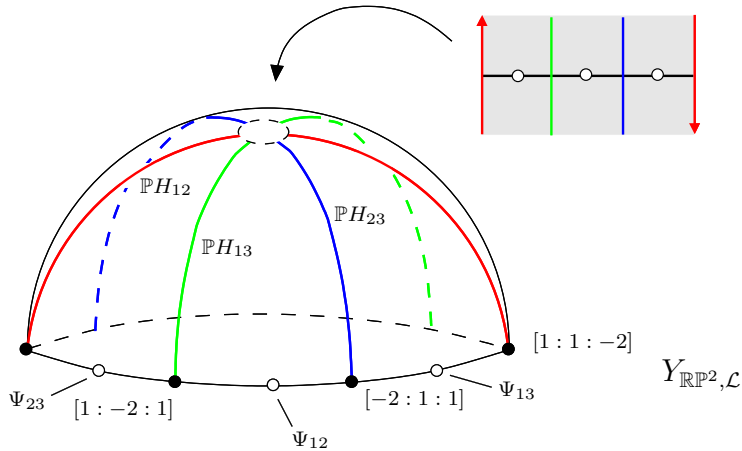

  \begin{picture}(0,0)%
    \includegraphics{blowup.pstex}%
  \end{picture}%
  \input{blowup.pstex_t}%
  
\caption{The wonderful model $Y_{\R\P^2, \LL}$.}
\label{fig_blowup}
\end{figure}

\end{explrm}

\begin{remrm}
  As already the low-dimensional Example~\ref{ex_S3P} shows, the
  $\LL$-stratification associated with the permutation action of
  $\Ss_n$ on $\R\P^{n{-}1}$ is different from the codimension~$2$
  truncation of the stabilizer stratification.
\end{remrm}



\end{document}